\documentclass[11pt]{article}
\usepackage{amsmath,amssymb,euscript,graphicx,amsfonts}
\usepackage{color}
\usepackage{enumerate}
\usepackage[dvips]{epsfig}
\usepackage{hyperref}

\newtheorem{theorem}{Theorem}[section]
\newtheorem{proposition}[theorem]{Proposition}
\newtheorem{lemma}[theorem]{Lemma}
\newtheorem{corollary}[theorem]{Corollary}

\newcommand{\Qed}{\rule{2.5mm}{3mm}}

\newcommand{\ZZ}{\mathbb{Z}}

\newcommand{\Ins}{\hbox{{\rm Ins}}}
\newcommand{\Out}{\hbox{{\rm Out}}}

\newenvironment{proof}{{\noindent \sc Proof.}}{\hfill $\Qed$ \\}

 \newcounter{case}
 
 \renewcommand{\thecase}{\arabic{case}}
 
 \newcounter{subcase}
 
\numberwithin{subcase}{case}

\def\ZZ{{\hbox{\sf Z\kern-.43emZ}}}

\input amssym.def

\input amssym.tex

\color{black}
\begin{document}


\begin{center}
{\bf\large ON CYCLIC  EDGE-CONNECTIVITY OF FULLERENES}
\end{center} 
\bigskip\noindent
\begin{center}

 Klavdija Kutnar{\small$^{a,}$}\footnotemark \, and  
 Dragan Maru\v si\v c{\small$^a$}$^,${\small$^b$}$^,$\addtocounter{footnote}{-1}\footnotemark$^,$*

\bigskip

{\it {\small  $^a$University of Primorska, Titov trg 4, 6000 Koper, Slovenia\\
$^b$University of Ljubljana, IMFM, Jadranska 19, 1000 Ljubljana, Slovenia\\}}
\
\end{center}

\addtocounter{footnote}{0}  \footnotetext{Supported in part by
``Agencija za raziskovalno dejavnost Republike Slovenije'', research program P1-0285.

~*Corresponding author e-mail: ~dragan.marusic@guest.arnes.si}

\begin{abstract}
A graph is  said to be {\em cyclic $k$-edge-connected}, 
if at least $k$ edges must be removed to disconnect 
it into two components,  each containing a cycle. 
Such a set of  $k$ edges 
is called a {\em cyclic-$k$-edge cutset}
and   it is called a
{\em trivial cyclic-$k$-edge cutset}  if at least one of 
the resulting two components induces a single $k$-cycle.

It is known that fullerenes, that is,  $3$-connected cubic planar graphs
all of whose faces are pentagons and hexagons, are cyclic  $5$-edge-connected.
In this article  it is shown that a fullerene $F$ containing a nontrivial cyclic-$5$-edge 
cutset admits two antipodal pentacaps, that is, two antipodal pentagonal 
faces whose neighboring faces are also pentagonal.
Moreover, it is shown that $F$
has a Hamilton cycle, and as a consequence 
at least $15\cdot 2^{\lfloor \frac{n}{20}\rfloor}$
perfect matchings, where $n$ is the order of $F$.
\end{abstract}

\bigskip

\begin{quotation}
\noindent {\em Keywords:} graph, fullerene graph, cyclic edge-connectivity, Hamilton cycle, perfect matching. 
\end{quotation}


\bigskip
\section{Introduction}
\label{sec:intro}
\indent

A {\em fullerene graph}  (in short a {\em fullerene}) 
is a  $3$-connected cubic planar graph,
all of whose faces are pentagons and hexagons. By Euler formula the number of pentagons equals $12$.
From a chemical point of view, fullerenes correspond to carbon 'sphere'-shaped 
molecules, the important class of molecules which is a basis of thousands of patents for a broad 
range of  commercial applications \cite{C2,C1}.
Graph-theoretic observations on structural properties of fullerenes are important in this
respect \cite{BM05,S02,S04,DDG98,JG1,JG2,JG3,JG4,AP01}.

In this paper cyclic  edge-connectivity of fullerenes, that is, the number $k$  such that
a fullerene cannot be separated into two components,
each containing a cycle, by deletion of fewer than $k$ edges, is considered. 
In general, a graph is  said to be {\em cyclically  $k$-edge-connected}, 
in short {\em cyclically} $k$-{\em connected},
if at least $k$ edges must be removed to disconnect 
it into two components,  each containing a cycle.  Cyclic 
edge-connectivity was extensively studied (see \cite{BJ91,K96,LH93,MS07,NS95}).
A set of  $k$ edges whose elimination disconnects a graph 
into two components, each containing a cycle,
is called a {\em cyclic-$k$-edge cutset}, in short a {\em cyclic-$k$-cutset}
and moreover, it is called a
{\em trivial cyclic-$k$-cutset} if at least one of 
the resulting two components induces a single $k$-cycle.
An edge from a cyclic-$k$-cutset is called {\em cyclic-cutedge}.

The concept of cyclic edge-connectivity played an important role in obtaining  some structural 
properties of fullerenes, such as bicriticality and $2$-extendability, that further
imply certain lower bounds on the number of perfect matchings  in fullerenes \cite{DD,D,D03, D06, ZZ}.
(Recall that a perfect matching $M$  in a graph  is a set of disjoint edges such that  
every vertex of the graph is covered by an edge from $M$.) 
A perfect matching in a graph  coincides, with the so called  Kekul\' e structure in chemistry,
and   the number of perfect matchings is an indicator of the stability of a fullerene.
However, in any fullerene  the cyclic edge-connectivity 
cannot exceed $5$, 
since by deleting the five edges connecting a pentagonal face,
two components each containing a cycle are obtained. 
In \cite[Theorem~2]{D03} it was proven that
it is in fact  precisely $5$. 

The main object of this paper is to give a more detailed 
description of cyclic-$5$-cutsets in fullerenes.
We show that, with the exception of a very special family of fullerenes possessing 
two antipodal pentacaps,
that is, two antipodal pentagonal faces whose neighboring
faces are also pentagonal, 
every other fullerene has only trivial cyclic-$5$-cutsets
(see Theorem~\ref{the:main}). 

Furthermore, we prove that fullerenes admitting a nontrivial cyclic-$5$-cutset  contain  a Hamilton
cycle, that is, a cycle going through all vertices  (see Theorem~\ref{the:main}), 
thus making a small contribution to the open problem regarding 
existence of Hamilton cycles in fullerenes \cite{M06}.
(Note that this problem is a special case of the Barnette's  conjecture~\cite{B69} which  says that 
every $3$-connected planar graph whose largest faces are hexagons,  contains a Hamilton cycle.) 
As an immediate consequence of this result it is shown that every fullerene of order $n$
admitting a nontrivial cyclic-$5$-cutset   has at least $15\cdot 2^{\lfloor \frac{n}{20}\rfloor}$
perfect matchings, thus improving (in the case of nontrivially cyclically $5$-edge-connected fullerenes) 
the best known lower bound  of $\lceil 3(m+2)/4 \rceil $ for the number of perfect matchings 
in a general fullerene of order $m$ (see  \cite{ZZ}).  




Throughout this paper graphs are finite, undirected
and connected,  unless specified otherwise.  
For notations and definitions not defined here we refer the reader to \cite{H69}.  
For adjacent vertices $u$ and $v$ in $X$,
we write $u \sim v$ and denote the corresponding edge by $uv$. 
Given a graph $X$ we let $V(X)$ and $E(X)$  
be the vertex set and the edge set  of $X$, respectively.
If $u\in V(X)$ then  $N(u)$ denotes the neighbors set of $u$ and
$N_i(u)$ denotes the set of vertices at distance  $i>1$  from $u$.
If $S\subseteq V(X)$, then $S^c=V(X)\setminus S$ denotes the complement of $S$ 
and the graph induced on $S$ is denoted by $X[S]$. Moreover, $X'=X[S+S']$
denotes the graph with the vertex set  $V(X')=S\cup S'$ and the edge set
$E(X')=E(X[S])\cup \{uv\mid u\in S, v\in S'\}$.


\section{Fullerenes admitting a nontrivial cyclic-$5$-cutset}
\label{sec:proof}
\indent

 An immediate consequence of  cyclic $5$-edge-connectivity of a fullerene
is the following result about the girth of a fullerene, that is, about the length of its 
smallest cycle.

\begin{proposition}
\label{pro:girth}
\label{cor:ful}  
The girth of a fullerene is $5$.
\end{proposition}

Let $C$ be a cycle in a  planar embedding of a fullerene $F$. 
Then we let the {\em inside} $\Ins(C)$ and the {\em outside} $\Out(C)$
be the set of vertices of $F$  that lie   inside $C$ and  outside $C$, respectively.
The following result is an immediate consequence of the well known Euler's formula 
for connected planar graphs. (Euler's formula states that the number of faces of a 
connected planar graph $X$ in its planar embedding 
is equal to $|E(X)|+|V(X)|-2$.)

\begin{proposition}
\label{pro:inside} Let $C$ be a cycle of length $10$ in a fullerene $F$ such that there exist exactly
five vertices on $C$ having a neighbor in the interior of $C$. 
Then exactly six pentagons exist   in the subgraph induced by $V(C) \cup \Ins(C)$. 
\end{proposition}

\begin{figure}[htbp]
\begin{footnotesize}
\begin{center}
\includegraphics[width=0.20\hsize]{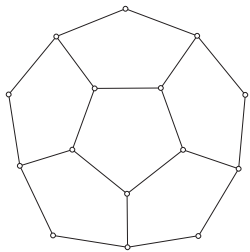}
\caption{\label{fig:cap}  The pentacap.}
\end{center}
\end{footnotesize}
\end{figure}

In this section we will prove that 
every fullerene admitting a nontrivial cyclic-$5$-cutset 
contains two antipodal pentacaps (see Theorem~\ref{the:main}),
where the {\em pentacap} is a planar graph on $15$ vertices 
with $7$ faces of which one  is a $10$-gon and six are pentagons 
(see Figure~\ref{fig:cap}). Observe
that the dodecahedron  is obtained as a union of two pentacaps, by
identifying the ten vertices on the outer ring  of the two pentacaps.
The following lemma  will be useful in this respect.

\begin{lemma}
\label{lem:c}
Let $F$ be a fullerene containing a ring $R$ of five faces, 
and let $C$ and $C'$ be the inner cycle and the outer cycle of $R$, respectively.
Then either 
\begin{enumerate}[(i)]
\itemsep=0pt
\item $C$ or $C'$ is a face, or 
\item both $C$ and $C'$ are of length $10$, and the five faces of $R$ are all hexagonal.
\end{enumerate}
\end{lemma}

\begin{proof}
Let  $f_0$, $f_1$, $f_2$, $f_3$ and $f_4$ be the faces in the ring $R$  
such that $f_i$ is adjacent to $f_{i+1}$, $i\in\ZZ_5$.
Let $T$ be the set of edges between these five faces
(depicted in bold in  Figure~\ref{fig:examples}). 
Clearly,  $T$ is a  cyclic-$5$-cutset of $F$. 
Moreover, if $T$ is a trivial cyclic-$5$-cutset of then  either $C$ or $C'$
is a face. We may therefore assume that $T$ is a nontrivial cyclic-$5$-cutset,
and that neither $C$ nor $C'$ is a face.

\begin{figure}[htbp]
\begin{footnotesize}
\begin{center}
\includegraphics[width=0.30\hsize]{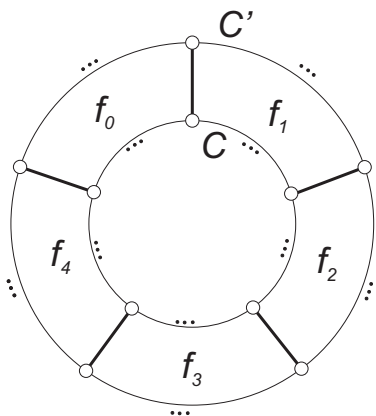}
\caption{{\small
\label{fig:examples}  The local structure of a fullerene that admits a nontrivial cyclic-$5$-cutset.}}
\end{center}
\end{footnotesize}
\end{figure}

Let $S = \Ins(C)$ and $S' = \Out(C')$, respectively,  be the 
 inside of $C$ and the outside of $C'$.
Let $l = l(C)$ and $l' = l(C')$  be the corresponding lengths of cycles $C$ and $C'$.
With no loss of generality, let $l \leq l'$.
Depending on whether the five faces are either all pentagonal at the one extreme or all hexagonal
at the other extreme, or possibly some pentagonal and some hexagonal, we have that
\begin{eqnarray}
\label{eq:sum}
15 \le l+ l'\le 20.
\end{eqnarray}
Since $T$ is a nontrivial cyclic-$5$-cutset and $C$ is not a face,
we must have that $l\geq 6$.

If $ l \in \{6,7\}$ then there either exists one edge or there are two edges having  one endvertex 
in  $C$ and the other in $S$  whose deletion disconnects  $F$, contradicting $3$-connectedness of $F$.
Therefore $l \geq 8$.

Suppose first that $l=8$. Then there exists a set of three edges 
$\{u_1v_1,u_2v_2,u_3v_3\}$ such that $u_1,u_2,u_3\in V(C)$, $u_i\ne u_j$ for $i\ne j$ 
and  $v_1,v_2,v_3\in S$. Since $F$ is 
cyclically $5$-connected it follows that the subgraph of $F$
induced  by  $S$ is a forest and so 
$X=F[S +\{u_1,u_2,u_3\}]$ is a forest, too.
Let $m=|S\cup\{u_1,u_2,u_3\}|$. 
Since all of the vertices in $S$ are of valency $3$ we have that
$2|E(X)| =  3(m-3) + 3$, whereas on the other hand $2|E(X)| = 2(|V(X)| -p) = 2(m-p)$, where 
$p$ denotes the number of components of $X$.  It follows that
\begin{eqnarray}
\label{eq:6}
m=6-2p.
\end{eqnarray}
Then by (\ref{eq:6}) it follows that $p=1$ and so  $m=4$. 
This implies that $|S|=1$ and therefore
$v_1=v_2=v_3$. But as $R$ consists of five faces 
one can easily see that  the subgraph induced by $V(C) \cup \Ins(C)$ 
contains a cycle of length $3$ or $4$,
contradicting Proposition~\ref{pro:girth}.
 
Suppose next that $l=9.$ Then there exists a set of four edges 
$\{u_1v_1,u_2v_2,u_3v_3,u_4v_4\}$ such that $u_1,u_2,u_3,u_4\in V(C)$, $u_i\ne u_j$ for $i\ne j$ 
and  $v_1,v_2,v_3,v_4\in S$. Since $F$ is 
cyclically $5$-connected it follows that the subgraph of $F$
induced  by $S$ is a forest and so 
$X=F[S +\{u_1,u_2,u_3, u_4\}]$ is a forest, too.
Let $m=|S\cup\{u_1,u_2,u_3,u_4\}|$. 
A counting argument similar to the one used in the previous paragraph 
gives us 
\begin{eqnarray}
\label{eq:8}
m=8-2p,
\end{eqnarray}
where  $p$ is the number of connected components of $X$. Clearly $p\le 2$. 
If $p=1$ then (\ref{eq:8}) gives us $m=6$ and so there exist two 
vertices in $S$, say  $v_1=v_2$  and $v_3=v_4$.  
Then $X$ has two vertices of 
valency $3$ and four vertices of valency $1$. 
But since  $l=9$ one can easily see that, as in the case $l=8$,
the subgraph induced by $V(C) \cup \Ins(C)$ 
contains a cycle of length $3$ or $4$,
contradicting Proposition~\ref{pro:girth}.
If $p=2$ then (\ref{eq:8}) implies that $m=4$, and so $S=\emptyset$. Then
without loss of generality $u_1\sim u_2$ and $u_3\sim u_4$.
But then again the fact that $l=9$ implies the existence of a cycle of length less then or equal to $4$
in the graph induced by   $V(C) \cup \Ins(C)$,
a contradiction.

Now (\ref{eq:sum}) implies that $l=l'=10$ and therefore all of the faces $f_i$, $i\in\ZZ_5$, on $R$
are hexagons, completing  the proof of Lemma~\ref{lem:c}.
\end{proof}


Given a ring of faces $R$ in a planar embedding of a fullerene $F$
with  inner cycle $C$ and   outer cycle $C'$ we let a face $f\in R$
be of type $(j)=(j)_C$ if there exist $j$ vertices on $f$ having the third neighbor  
(the one different form the immediate neighbors on $C$) in
$\Ins(C)\cup V(C)$. Clearly $0\le j\le 2$. Further, a ring $R$
of $r$ faces $f_1,f_2,\ldots,f_r$ such that $f_{i}$ is adjacent to $f_{i+1}$, $i\in\ZZ_r$,
is said to be of type $(j_1\,j_2\,\ldots\, j_r)=(j_1\,j_2\,\ldots\, j_r)_C$ if  $f_i\in R$ is of type
$(j_i)$ for every $i\in\ZZ_r$. For example, a ring in which the inner cycle is a face
is of type $(00000)$ if the inner cycle is a pentagon and 
of type $(000000)$ if the inner cycle is a hexagon.


We may now prove the main theorem of this paper.

\medskip

\begin{theorem}
\label{the:main} 
Let $F$ be a fullerene admitting a nontrivial cyclic-$5$-cutset.
Then $F$ contains a pentacap, more precisely, 
either $F$ is the dodecahedron or it contains two disjoint antipodal  pentacaps.
\end{theorem}

\begin{proof}
It is clear that the dodecahedron admits a nontrivial cyclic-$5$-cutset. 
Therefore, let $F$ be a fullerene  admitting  a nontrivial cyclic-$5$-cutset $T$ different from 
the dodecahedron. Then there exists a ring $R$ of five faces $f_0$, $f_1$, $f_2$, $f_3$ and $f_4$ in $F$
such that $f_i$ is adjacent to $f_{i+1}$, $i\in\ZZ_5$, via an edge from $T$.
Let $C_0$ be the inner cycle of $R$ and $C_1$ be the outer cycle of $R$. 
Let $S_0=\Ins(C_0)$ and $S_1=\Out(C_1)$, respectively,  be the 
inside of $C_0$ and the outside of $C_1$.
Let $l_0=l(C)$ and $l_1= l(C')$  be the corresponding lengths of cycles $C_0$ and $C_1$.
By Lemma~\ref{lem:c} we have that $l_0=l_1=10$.   
Hence the five faces of $R$ are all hexagonal and  precisely five vertices  
having a neighbor in $S_0$ exist on $C_0$. Depending on the arrangement of 
the vertices on $C_0$ having a neighbor in $S_0$ the ring $R$ is of one of the following six types:
$(01112)$, $(01121)$, $(00212)$, $(00122)$, $(02102)$ or $(11111)$ 
(see Figure~\ref{fig:cases}). (Note that these are  
all possible types  since $l_0=10$.)  

\begin{figure}[htbp]
\begin{footnotesize}
\begin{center}
\includegraphics[width=0.90\hsize]{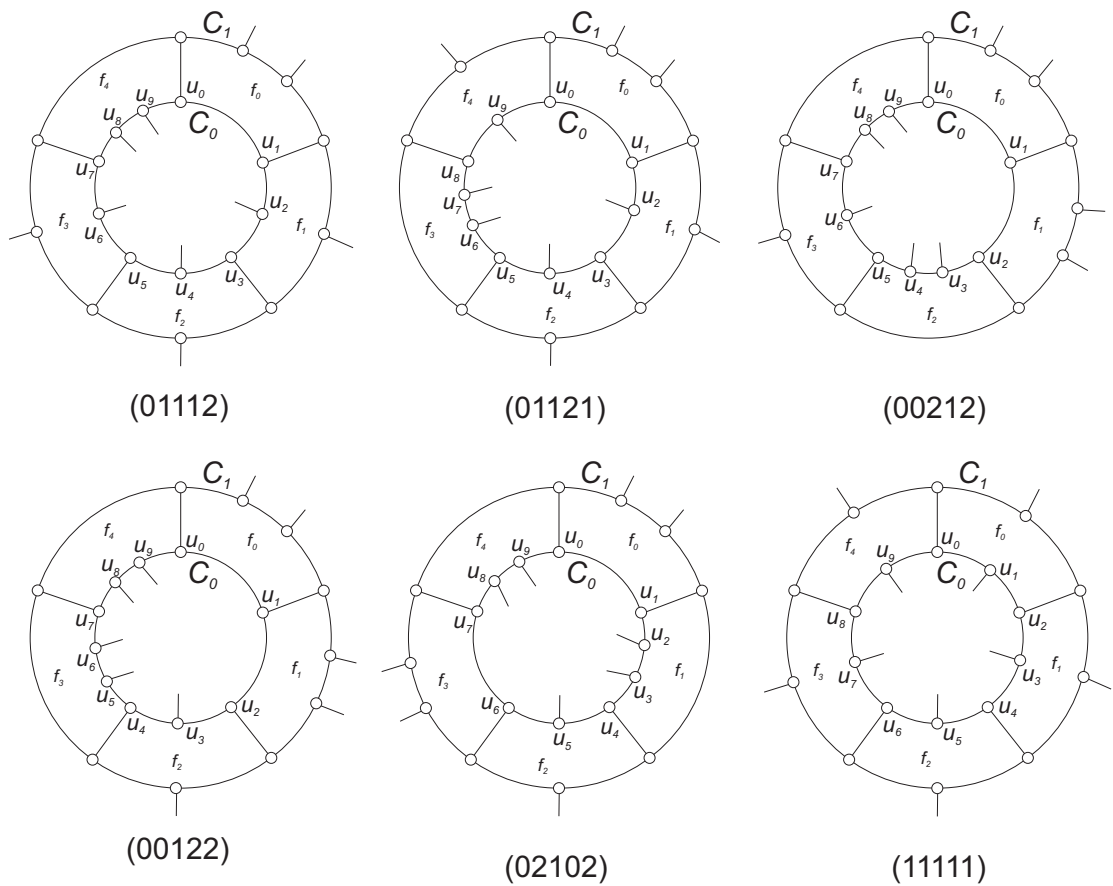}
\caption{{\small\label{fig:cases}  The six possible types of   $R$.}}
\end{center}
\end{footnotesize}
\end{figure}

We claim that only  type (11111)  can occur.
Let $C_0=u_0u_1\ldots u_9$ and first, suppose 
that there exist $i,j\in\ZZ_9$, $j\not\in\{i-1,i,i+1\}$, such that $u_i$ and $u_j$ 
are adjacent. Since  $F$ has   girth  $5$ we have that
$$
4\le i-j \le 5.
$$
Now planarity and $3$-connectivity of $F$  combined together imply that
either $u_i,u_{i+1},\ldots,u_{j-1},u_j$ or $u_j,u_{j+1},\ldots,u_{i-1},u_i$ determinate
a face, say $f$, in $F$ (observe that $u_i$ and $u_j$ 
lie on nonneighboring faces of $R$). 
Checking all possible types one can see that  
$R$ is either of type $(00212)$ or of type $(00122)$. But then 
$ff_if_{i+1}f_{i+2}$   (where $f_i$, $f_{i+1}$ and $f_{i+2}$, $i\in \ZZ_5$,
are   faces of $R$ of nonzero type), is a ring of four faces. 
This is impossible in view of cyclic $5$-edge-connectivity of $F$. 

Next, suppose that there exist $i,j\in\ZZ_9$, $j\ne i$, such that $N(u_i)\cap N(u_j)\cap S_0\ne \emptyset$.  
Since  $F$ has   girth  $5$ we have that
$$
3\le i-j \le 5.
$$
Now planarity and $3$-connectivity of $F$  combined together imply that
either $u_i,u_{i+1},\ldots,u_{j-1},u_j$ or $u_j,u_{j+1},\ldots,u_{i-1},u_i$ determinate
a face, say $f$, in $F$, and so $u_i$ and $u_j$ lie on nonneighboring faces of $R$.
Clearly, $f$ is a neighbor of either one or two faces of $R$ having type $(0)$.
If the former case holds then $f$ is pentagonal, and moreover $f$ together with the four faces in $R$
different from the type $(0)$ face adjacent to $f$, form a ring of five faces, contradicting
Lemma~\ref{lem:c}. If the latter holds then $f$ together with the three nonzero type ($\ne (0)$) faces 
of $R$ form a ring of four faces, contradicting cyclic  $5$-edge-connectivity of $F$.

Hence we have that $N(u_i)\cap N(u_j)\cap S_0= \emptyset$ for $i,j\in\ZZ_9$, $j\ne i$,  
and therefore $|S_0|\ge 5$.  Note also that this implies that $R$ is of type different
from type $(00212)$ and type $(00122)$ (see also Figure~\ref{fig:cases}).
%
%
If  $F[S_0]$ and so also $X=F[S_0+\{u_1,u_2,u_3,u_4,u_5\}]$  is a forest then
a simple counting argument shows that 
\begin{eqnarray}
\label{eq:10}
|S_0\cup\{u_1,u_2,u_3,u_4,u_5\}|=10-2p,
\end{eqnarray}
where  $p$ is the number of connected components of $X$. 
But since $|S_0\cup\{u_1,u_2,u_3,u_4,u_5\}|\ge 10$ and $p\ge 1$ this is clearly impossible.
Therefore  $F[S_0]$ contains a cycle and  
there exists a ring $R'$ of five faces whose outer cycle
is $C_0$. By Lemma~\ref{lem:c} either the inner cycle of $R'$ is 
a pentagonal face and all faces of $R'$
are pentagonal, or all faces of $R'$ are hexagonal.
In the first case we are done.
In the second case we replace $R$ with $R'$ in our analyze.
Continuing with this line of argument and using the fact that $F$
is finite in some stage we have to reach a ring of five pentagonal faces
giving rise to the pentacap. 
In particular, since a ring adjacent to a ring of type  $(01112)$, $(01121)$ or $(02102)$ 
always contains a hexagonal face this process can stop only if $R$ is of type $(11111)$.
This completes the proof of Theorem~\ref{the:main}.
\end{proof}

We remark that Theorem~\ref{the:main} implies that fullerenes 
admitting a nontrivial cyclic-$5$-cutset are a special class of the so-called carbon nanotubes
(see \cite{H99}).


\section{Hamilton cycles in fullerenes admitting a nontrivial cyclic-$5$-cutset}
\label{sec:ham}

\indent

In this section  it is proven that the Barnette conjecture 
is true for fullerenes admitting a nontrivial cyclic-$5$-cutset (see Theorem~\ref{the:ham}).
The key factor in the proof of this result is Theorem~\ref{the:main}. In particular,
by Theorem~\ref{the:main} we know that each fullerene admitting a nontrivial cyclic-$5$-cutset
contains two antipodal pentacaps and moreover, from the proof of Theorem~\ref{the:main} one may
deduct that between these two pentacaps there exist rings of five hexagonal faces 
(hexagonal rings) such that in each hexagon $H$ exist two vertices, first of which 
has a neighbor inside the ring which $H$ belongs to
and second have a neighbor outside the ring which $H$ belongs to.

\begin{theorem}
\label{the:ham}
Let $F$ be a fullerene admitting a nontrivial cyclic-$5$-cutset. 
Then $F$ has a Hamilton cycle. Moreover, 
\begin{enumerate}[(i)]
\itemsep=0pt
\item if the number of hexagonal faces in $F$ is 
odd then there exists a path of faces containing precisely two pentagons from each
of the two pentacaps in $F$ whose boundary gives rise to a Hamilton cycle in $F$; and 
\item if the number of hexagonal faces in $F$ is 
even then there exists a path of faces containing precisely six pentagons, of which two
are from the first pentacap and four from the other pentacap,
whose boundary gives rise to a Hamilton cycle in $F$.
\end{enumerate}
\end{theorem}

\begin{figure}[h!]
\begin{footnotesize}
\begin{center}
\includegraphics[width=0.700\hsize]{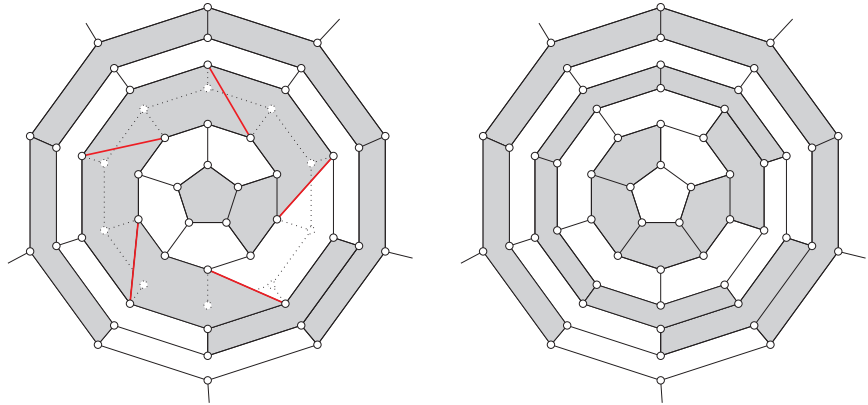}
\caption{{\small \label{fig:even}  A local structure of a Hamilton cycle in $\bar{F}$ on the 
left-hand side  picture and a local structure of a Hamilton cycle in $F$ on the
right-hand side  picture for $r+1$  even.}}
\end{center}
\end{footnotesize}
\end{figure}

\begin{figure}[h!]
\begin{footnotesize}
\begin{center}
\includegraphics[width=0.700\hsize]{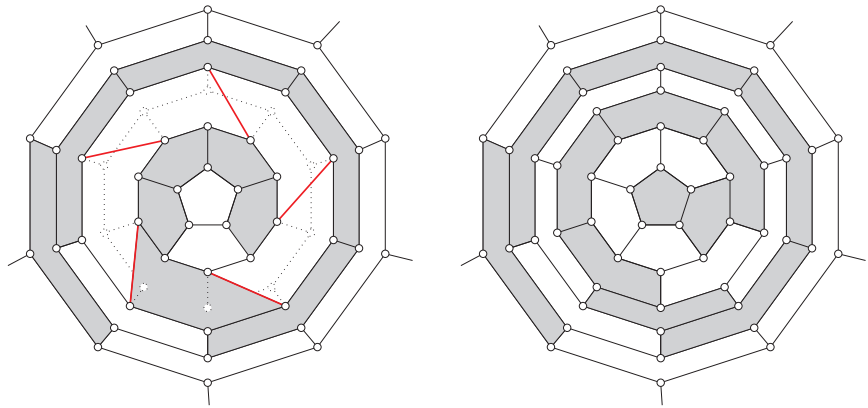}
\caption{{\small\label{fig:odd}  A local structure of a Hamilton cycle in $\bar{F}$ on the 
left-hand side  picture and  a local structure of a  Hamilton cycle in $F$ on the 
right-hand side  picture for $r+1$  odd.}}
\end{center}
\end{footnotesize}
\end{figure}

\begin{proof}
By Theorem~\ref{the:main} the fullerene $F$ contains two antipodal pentacaps.
From the proof of Theorem~\ref{the:main}  we may deduct that between
these two pentacaps there exist  rings each consisting of five hexagonal faces, in short
hexagonal rings, such that 
for each of these rings the following holds: 
in each hexagon in the ring there exist a vertex having a neighbor inside this ring, and
a vertex having a neighbor outside this ring.

Let $k$ be the number of hexagonal faces in $F$. Observe that $k=5r$ where
$r$ is the number of hexagonal rings in $F$. We proceed
 the proof by induction on $r$.

If $r=0$ then $F$ is the dodecahedron in which a path of faces containing 
precisely six pentagons whose boundary gives rise to a Hamilton cycle clearly exists
(see also Figure~\ref{fig:dod}). Further, in Figures~\ref{fig:1} and~\ref{fig:2}
Hamilton cycles in $F$ are shown if $r=1$ and $r=2$. 
Hence, the statement of the theorem holds  for $r\le 2$.

\begin{figure}[h!]
\begin{footnotesize}
\begin{center}
\includegraphics[width=0.20\hsize]{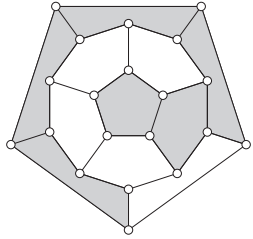}
\caption{{\small\label{fig:dod}  A Hamilton cycle in the dodecahedron ($r=0$).}}
\end{center}
\end{footnotesize}
\end{figure}

\begin{figure}[htb]
\begin{minipage}[b]{0.5\linewidth}
\centering
\includegraphics[width=0.500\hsize]{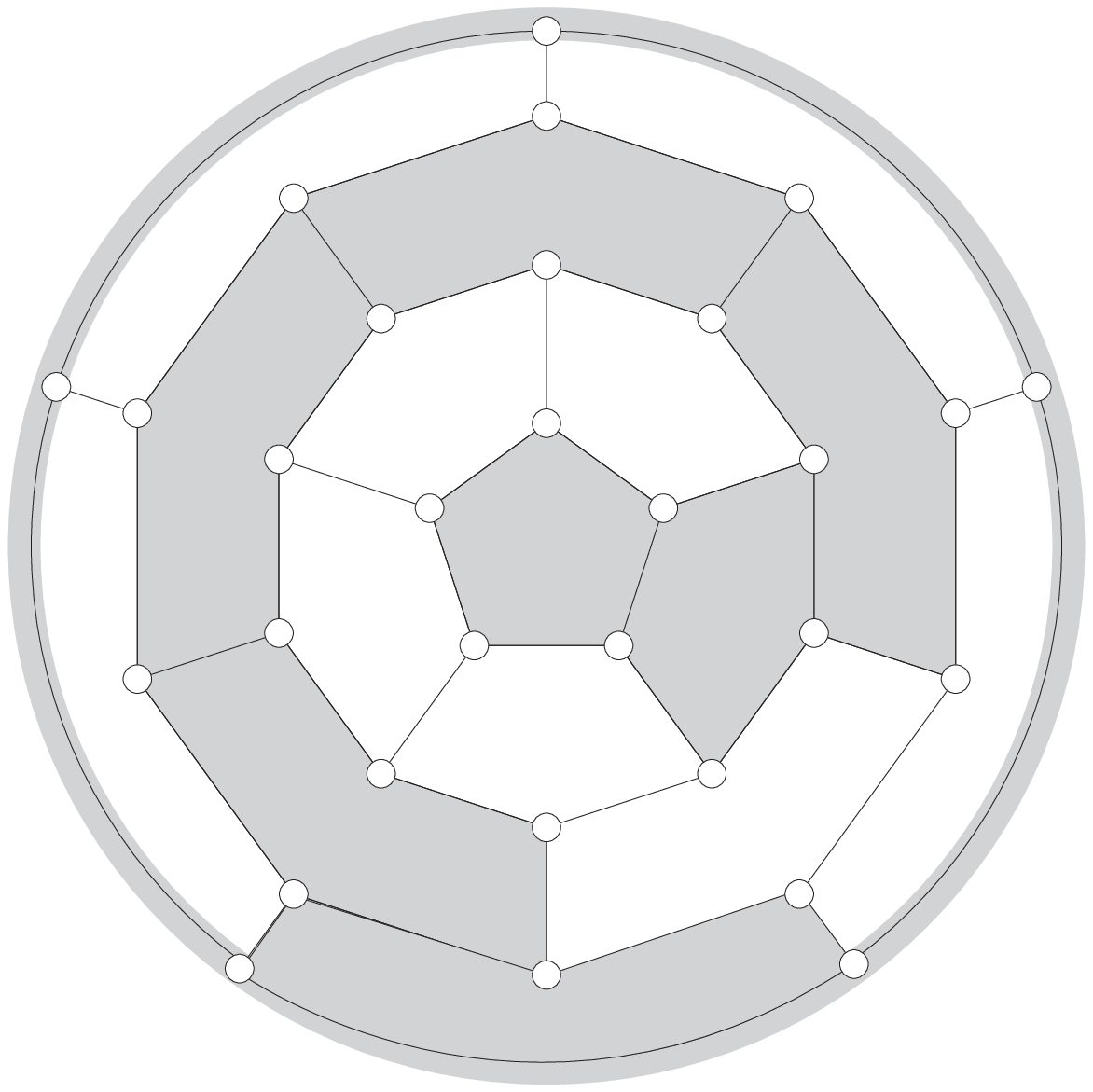}
\caption{{\small\label{fig:1}  A Hamilton cycle in $F$ if $r=1$.}}
\end{minipage}\hspace{0.5cm}
\begin{minipage}[b]{0.5\linewidth}
\centering
\includegraphics[width=0.550\hsize]{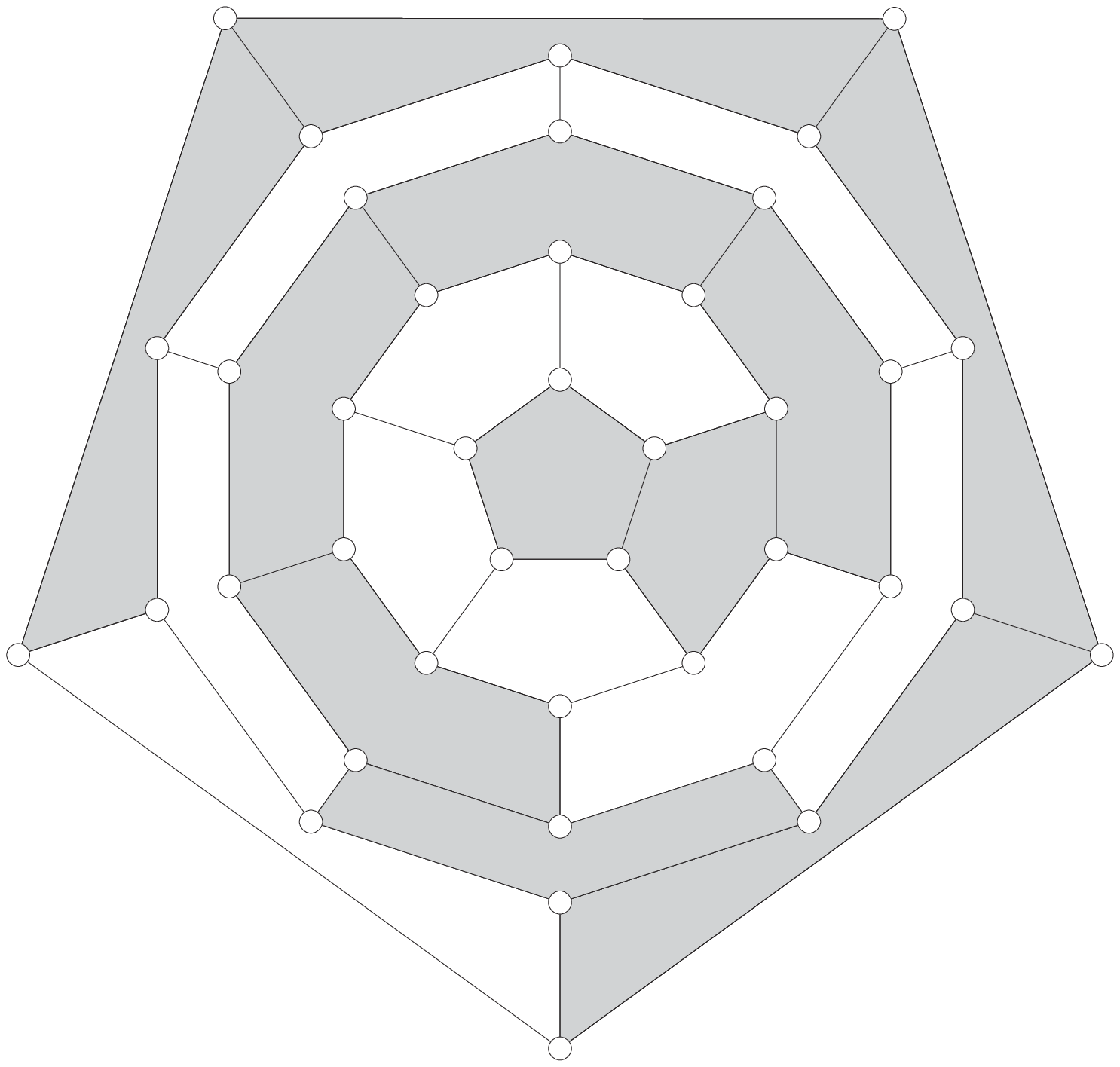}
\caption{{\small\label{fig:2}  A Hamilton cycle in $F$ if  $r=2$.}}
\end{minipage}
\end{figure}

Assume now that the statement of the theorem holds for fullerenes with $r>2$ hexagonal rings and
let $F$ be a fullerene admitting a nontrivial cyclic-$5$-cutset with $r+1$ hexagonal rings.
Denote by $R$ the hexagonal ring in $F$ adjacent to one of the two
pentacaps, and denote by $R'$ the hexagonal ring adjacent to $R$. Furthermore,
let $C_0$ and $C_1$ be the inner and the outer cycle of $R$, respectively.
Clearly $C_1$ is the inner cycle of $R'$. Let $C_2$ be the outer cycle of $R'$. 
Further, let $v_i^j\in V(C_i)$, $i\in \{0,1,2\}$ and $j\in \ZZ_{10}$, be such that
$v_0^{2j}\sim v_1^{2j}$ and $v_1^{2j+1}\sim v_2^{2j+1}$. Now we construct
a fullerene admitting a nontrivial cyclic-$5$-cutset with $r$ hexagonal rings as follows.
By deleting of all the vertices on the cycle $C_1$, 
that is vertices $v_1^j$, $j\in\ZZ_{10}$ (with the incident edges),
and by adding edges $v_0^{2j}v_2^{2j-1}$, $j\in\ZZ_{10}$, we obtain a fullerene
admitting a nontrivial cyclic-$5$-cutset with $r$ hexagonal rings.
Denote this fullerene  by $\bar{F}$ (see Figures~\ref{fig:even} and~\ref{fig:odd}).  

If $r+1$ is even then $r$ is odd.  Hence, by the induction hypothesis a path of 
faces containing precisely two pentagons from each
of the two pentacaps  whose boundary gives rise to a Hamilton cycle exists in $\bar{F}$.
But then one can  construct a path of faces in $F$ containing precisely six pentagons 
whose boundary gives rise to a Hamilton cycle, as illustrated in 
Figure~\ref{fig:even}.

Suppose now that $r+1$ odd, $\bar{F}$ has an even number of hexagons.
Therefore, by induction hypothesis there exists
a path of faces containing precisely six pentagons, of which two
are from the first pentacap and four from the other pentacap,
whose boundary gives rise to a Hamilton cycle in $\bar{F}$.
Again, one can  construct a path of 
faces in $F$ containing precisely two pentagons from each of the two pentacaps 
whose boundary gives rise to a Hamilton cycle  in $F$ as illustrated in Figure~\ref{fig:odd}.
This completes the proof Theorem~\ref{the:ham}.
\end{proof}

\begin{figure}[h!]
\centering
\includegraphics[width=0.590\hsize]{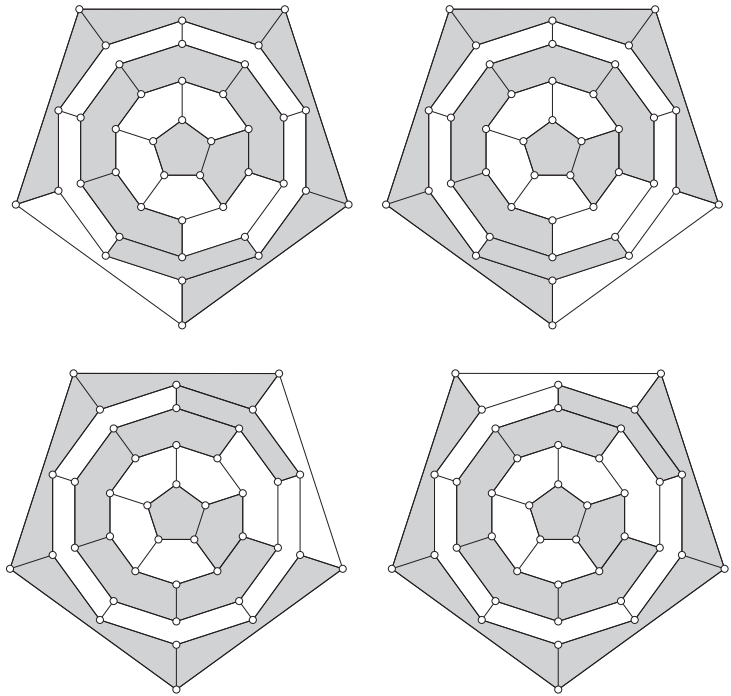}
\caption{{\small\label{fig:2iz}  Four different paths of faces in a fullerene $F$
admitting a nontrivial cyclic-$5$-cutset with $r=2$ whose boundary gives rise
to a Hamitlon cycle in $F$ for a choosen pentagon adjacent
to the central pentagon in a pentacap. 
Since the pentagon adjacent
to the central pentagon in the pentacap can be choosen in five differene ways this 
fullerene has at least $5\cdot 2^{\frac{r}{2}+1}=20$ different Hamilton cycles.}}
\end{figure}

Observe that a path of faces in a fullerene admitting a nontrivial cyclic-$5$-cutset $F$
whose boundary  gives rise to a Hamilton cycle in $F$ that was constructed in the proof
of Theorem~\ref{the:ham} is not unique.
As illustrated in Figure~\ref{fig:2iz}
one can  see that the following proposition holds.

\begin{proposition}
\label{pro:hx}
Let $F$ be a fullerene of order $n$ admitting a nontrivial cyclic-$5$-cutset
and let $r$ be the number of hexagonal rings in $F$. 
Then  
\begin{enumerate}[(i)]
\itemsep=0pt
\item if $r$ is even then $F$ has at least $5\cdot 2^{\frac{r}{2}+1}$ different Hamilton cycles; and
\item if $r$ is odd then $F$ has at least $5\cdot 2^{\frac{r+1}{2}}$ different Hamilton cycles.
\end{enumerate}
\end{proposition}

Since  every Hamilton cycle in a fullerene   $F$ gives rise
to three perfect matchings in $F$, 
Proposition~\ref{pro:hx} gives a  lower bound of the number on perfect 
matchings in a  fullerene admitting a nontrivial cyclic-$5$-cutset.
Since the number of hexagonal rings in a fullerene addmiting a nontrivial cyclic-$5$-cutset
is $r=\frac{n}{10}-2$ the following corollary holds.
 
\begin{corollary}
\label{cor:pm}
Let $F$ be a fullerene of order $n$ admitting a nontrivial cyclic-$5$-cutset. Then
the number of perfect 
matchings in $F$ is at least $15\cdot 2^{\lfloor \frac{n}{20}\rfloor}$. 
\end{corollary}


\end{document}